\documentclass[a4paper,12pt]{amsart}

\usepackage{float}
\usepackage{euscript,eufrak,verbatim}
\usepackage{graphicx}
\usepackage[usenames]{color}
\usepackage[colorlinks,linkcolor=red,anchorcolor=blue,citecolor=blue]{hyperref}
\usepackage{amsmath}
\usepackage{amsthm}
\usepackage[all]{xy}
\usepackage{graphicx} 
\usepackage{amssymb}

\theoremstyle{plain}
\newtheorem{main theorem}{Main Theorem}
\newtheorem{theorem}{Theorem}[section]

\theoremstyle{definition}

\makeatletter 
\@addtoreset{equation}{section}
\numberwithin{equation}{section}

\newcommand{\norm}[1]{\left\lVert#1\right\rVert}

\title{Simple proof of the global inverse function theorem via the Hopf--Rinow theorem}

\author{Shinobu Ohkita}

\address{Department of Mathematics, Kyoto University, Kitashirakawa Oiwake-cho, Sakyo-ku, Kyoto 606-8502, Japan}

\email{okita.shinobu.42w@st.kyoto-u.ac.jp}

\author{Masaki Tsukamoto}

\address
{Department of Mathematics, Kyoto University, Kitashirakawa Oiwake-cho, Sakyo-ku, Kyoto 606-8502, Japan}

\email{tsukamoto@math.kyoto-u.ac.jp}

\begin{document}



\thanks{M.T. was supported by JSPS KAKENHI JP21K03227.}

\maketitle

\begin{abstract}
We explain that Hadamard’s global inverse function theorem very simply follows from 
the Hopf--Rinow theorem in Riemannian geometry.
\end{abstract}

\section{Global inverse function theorem}  \label{section: global inverse function theorem}

The purpose of this note is to explain a very simple proof of Hadamard’s global inverse function theorem.
For a smooth map $f = (f_1, \dots, f_n)\colon \mathbb{R}^n\to \mathbb{R}^n$ and $x\in \mathbb{R}^n$ 
we denote by $Df(x) = \left(\partial f_i/\partial x_j\right)_{i,j}$ the Jacobian matrix of $f$ at $x$.
For an $n\times n$ matrix $A$ we denote its operator norm by $\norm{A}$.

The following statement is known as the \textbf{global inverse function theorem}:

\begin{theorem}[Global inverse function theorem] \label{theorem: global inverse function theorem}
Let $f\colon \mathbb{R}^n\to \mathbb{R}^n$ be a smooth map such that the Jacobian matrix $Df(x)$ is invertible 
at every $x\in \mathbb{R}^n$ and that $\norm{\left(Df(x)\right)^{-1}}$ is bounded over $x\in \mathbb{R}^n$.
Then $f$ is a diffeomorphism of $\mathbb{R}^n$ onto itself.
\end{theorem}

This theorem is usually attributed to the paper of Hadamard \cite{Hadamard} in the early 20th century.
Several authors \cite{Schwartz, Gordon, Plastock, Miller, Rabier, RS15} 
had provided its proofs and generalizations from various viewpoints. 
The global inverse function theorem was used for the study of 
oscillatory integral in \cite[p. 302]{Asada--Fujiwara}.
We can prove the Fundamental Theorem of Algebra by using the global inverse function theorem; see
the book of Krantz--Parks \cite[\S 6.2, Theorem 6.2.5]{Krantz--Parks}.

The point of Theorem \ref{theorem: global inverse function theorem}
is the assumption that $\norm{\left(Df(x)\right)^{-1}}$ is bounded over $\mathbb{R}^n$.
If we remove it, then the statement is false in general. For example, consider the (complex) exponential function
$\mathbb{C}\ni z\mapsto e^x\in \mathbb{C}$. Its Jacobian matrix (as a map from $\mathbb{R}^2$ to $\mathbb{R}^2$) is 
\begin{equation} \label{eq: Jacobian matrix}
     e^x  \begin{pmatrix} \cos y & -\sin y \\  \sin y  &  \cos y \end{pmatrix},  \quad (z = x + y\sqrt{-1}) 
\end{equation}     
which is invertible everywhere.  But $e^z$ is neither injective nor surjective.
Indeed the inverse of the matrix (\ref{eq: Jacobian matrix}) is given by
\[   e^{-x} \begin{pmatrix} \cos y & \sin y \\  -\sin y & \cos y\end{pmatrix}. \]
Its operator norm is $e^{-x}$, which is unbounded over $\mathbb{C}$.

We will show that Theorem \ref{theorem: global inverse function theorem} 
is a very easy consequence of the \textbf{Hopf--Rinow theorem}.
The Hopf--Rinow theorem is a standard result of Riemannian geometry and (a part of) its statement is as follows (see e.g. \cite[Theorems 4.1 and 4.2]{Kobayashi--Nomizu}):

\begin{theorem}[Hopf--Rinow theorem]
Let $M$ be a connected Riemannain manifold, and suppose that it is complete as a metric space 
(i.e., every Cauchy sequence converges).
Then we have:
  \begin{enumerate}
    \item  For any two points in $M$ there exists a length-minimizing geodesic between them.
    \item  For every point $p\in M$ the exponential map $\exp_p$ is defined all over the tangent space $T_p M$.
  \end{enumerate}
\end{theorem}

\begin{proof}[Proof of Theorem \ref{theorem: global inverse function theorem}]
It is enough to show that $f$ is a bijection from $\mathbb{R}^n$ onto $\mathbb{R}^n$.
Let $h$ be the Euclidean metric on $\mathbb{R}^n$.
From the assumption, there is $c>0$ such that $h\left(df(v), df(v)\right)\geq c\cdot h(v, v)$ 
for all tangent vectors $v\in T\mathbb{R}^n$.
We define a Riemannian metric $g$ on $\mathbb{R}^n$ by 
$g(u, v) = h\left(df_x(u), df_x(v)\right)$ for $u, v\in T_x \mathbb{R}^n$.
The map $f\colon (\mathbb{R}^n, g) \to (\mathbb{R}^n, h)$ is a local isometry (i.e. the pull-back of $h$ by $f$ is equal to $g$).

Let $d$ be a distance function induced by $g$.
Then we have $d(x, y) \geq \sqrt{c}\, |x-y|$, where $|x-y|$ is the standard Euclidean distance.
In particular $(\mathbb{R}^n, d)$ is complete as a metric space. So we can apply the Hopf--Rinow theorem
to $(\mathbb{R}^n, g)$.

\begin{itemize}
  \item We prove that $f$ is injective. Suppose $x$ and $y$ are two distinct points of $\mathbb{R}^n$.
By the Hopf--Rinow theorem there exists a geodesic $\gamma$ between them (with respect to the metric $g$).
Since $f$ is a local isometry, 
$f\circ \gamma$ is a geodesic in $(\mathbb{R}^n, h)$. 
So it is the line segment between $f(x)$ and $f(y)$.
The image of $f\circ \gamma$ is not a single point because $df_x\colon T_x\mathbb{R}^n\to T_{f(x)}\mathbb{R}^n$ is invertible. 
Therefore $f(x) \neq f(y)$.

   \item  Next we prove that $f$ is surjective. We assume $f(0) = 0$ for simplicity.
   Let $x\in \mathbb{R}^n$ be an arbitrary point. We define a line $\ell\colon \mathbb{R}\to \mathbb{R}^n$
   by $\ell(t) = tx$. 
   Let $v:= \ell^\prime(0)$ be the tangent vector of $\ell$ at the origin.
   Take $u\in T_0\mathbb{R}^n$ with $df_0(u) = v$, and consider a geodesic 
   $\gamma(t) := \exp_0(t u)$ (with respect to the metric $g$).
   This is defined for all $t\in \mathbb{R}$ by the Hopf--Rinow theorem.
   Then $f\circ \gamma$ is a line of constant speed and $\left(f\circ \gamma\right)^\prime(0) = v = \ell^\prime(0)$.
   Thus $f\circ \gamma(t) = \ell(t)$. In particular we have $f\left(\gamma(1)\right) = x$.
\end{itemize}
This finishes the proof of Theorem \ref{theorem: global inverse function theorem}.
\end{proof}

\section{Discussions and remarks}  \label{section: discussions and remarks}

\begin{enumerate}
   \item The above proof of the surjectivity of $f$ is close to the standard proof of the following 
   theorem of Riemannian geometry \cite[Theorem 4.6]{Kobayashi--Nomizu}: \textit{Let $M$ and $N$ be connected Riemannian 
   manifolds of the same dimension, and suppose $M$ is complete. Then any local isometry from $M$ to $N$ is a covering map.}
   (See also the book of Hermann \cite[Chapter 21]{Hermann}.)
   Therefore it is fair to say that Riemannian geometry provides a natural (and, in a sense, broader) framework for the global inverse function theorem.
   More general global inverse function theorems were developed by Guti\'{e}rrez--Biasi--Santos \cite{GBS09}
   in the context of Riemannian geometry. 
   Rabier \cite{Rabier} developed a global inverse function theorem for Finsler manifolds.

   \item As far as the authors know, most known proofs of the global inverse function theorem use some topological 
   argument (e.g. the covering space theory in \cite{Plastock} or homotopy arguments in \cite{Schwartz, Miller}). 
   But the above proof does not use any topological 
   argument.
   Moreover, the above proof of the surjectivity of $f$ is constructive, at least in principle.
   If the map $f$ is given by an explicit formula,
   then the geodesic $\gamma(t)$ is a solution of an explicit ordinary differential equation.
   So we can numerically solve it. In particular we can approximately calculate the point $\gamma(1)$ satisfying
   $f\left(\gamma(1)\right) = x$.

   \item  
   The metric $g$ was defined by $g(u, v) = h\left(df(u), df(v)\right)$ in the above proof.
   The geodesic equation involves first order differentials of the metric $g$ and hence second order differentials of 
   the map $f$. So we need to assume that $f$ is at least $C^2$ maps.
   However it is known that the global inverse function theorem holds for $C^1$ maps $f$ as well \cite[p.16]{Schwartz}.
   Therefore our proof does not provide an optimal regularity.
   But we think that the striking simplicity of the proof compensates for this drawback.

   \item 
   The proof of the Hopf--Rinow theorem is not very easy.
   So one can argue that our proof just translates one difficulty into another.
   (In some sense, our discovery is that the Hopf--Rinow theorem contains all the ingredients needed for the proof of 
   the global inverse function theorem.)
   But the Hopf--Rinow theorem is certainly much more well-known than the global inverse function theorem.
   We think that it is nice to see that the Hopf--Rinow theorem has such an unexpected application.

   \item  
   Plastock \cite[\S 3]{Plastock} studied generalizations of the global inverse function theorem motivated by the Hopf--Rinow theorem.
   So our paper is conceptually similar to \cite{Plastock}. 
   But our viewpoint is different from \cite{Plastock}.
   The point of our paper is that the global inverse function theorem is a corollary of the Hopf--Rinow theorem.
   On the other hand, Plastock \cite[\S 3]{Plastock} did not use the Hopf--Rinow theorem itself, but he developed a method motivated by   
   the proof of the Hopf--Rinow theorem.

   \item Probably the most useful “application” of this paper is to use it for education.
   It may be nice to explain the global inverse function theorem as a corollary of the Hopf--Rinow theorem 
   in an introductory course on Riemannian geometry.
   Indeed the content of this paper grew out of an undergraduate course conducted by the second named author in Kyoto university.
   He posed a problem concerning the global inverse function theorem in an exercise course on 
   geometry. Then the first named author found a proof using the Hopf--Rinow theorem.
   His idea looked so beautiful that the second named author thought that it should be published.
   Therefore all the ideas of this paper are due to the first named author.
   The second named author just elaborated technical details and presentations.
\end{enumerate}


\begin{thebibliography}{99}

\bibitem[AF78]{Asada--Fujiwara}
K.~Asada, D.~Fujiwara,
On some oscillatory integral transformations in $L^2(\mathbb{R}^n)$,
 Jpn. J. Math.  \textbf{4} (1978) 299–361.


\bibitem[Gor72]{Gordon}
W.~B.~Gordon,
On the diffeomorphisms of Euclidean space,
Amer. Math. Monthly \textbf{79} (1972) 755-759.


\bibitem[GBS09]{GBS09}
C.~Guti\'{e}rrez, C.~ Biasi, E.~L.~d.~Santos,
Global inverse mapping theorems,
Candernos De Mathem\'{a}tica \textbf{10} (2009) 9-18, 
Artigo N\'{u}mero SMA\#311.



\bibitem[Had04]{Hadamard}
J.~Hadamard, 
Sur les transformations ponctuelles,
Bull. Soc. Math. France \textbf{34} (1904) 71-84.



\bibitem[Her77]{Hermann}
R.~Hermann,
Differential geometry and the calculus of variations,
Second edition, Interdisciplinary Mathematics, Vol. XVII, Math Sci Press, Brookline, Mass. (1977)




\bibitem[KN63]{Kobayashi--Nomizu}
S.~Kobayashi, K.~Nomizu,
Foundations of differential geometry, Volume I, 
John Wiley \& Sons, Inc.,1963.


\bibitem[KP02]{Krantz--Parks}
S.~G.~Krantz, H.~R.~Parks,
The implicit function theorem, History, theory, and applications,
Birkh\"{a}user Boston (2002)



\bibitem[Mil84]{Miller}
J.~D.~Miller, 
Some global inverse function theorems,
J. Math. Anal. Appl. \textbf{100} (1984) 375-384.



\bibitem[Pla74]{Plastock}
R.~Plastock,
Homeomorphisms between Banach spaces,
Trans. Amer. Math. Soc. \textbf{200} (1974) 169-183.



\bibitem[Rab97]{Rabier}
P.~J.~Rabier,
Ehresmann fibrations and Palais--Smale conditions for 
morphisms of Finsler manifolds, 
Ann. of Math. \textbf{146} (1997) 647-691.




\bibitem[RS15]{RS15}
M.~Ruzhansky, M.~Sugimoto,
On global inversion of homogeneous maps,
Bull. Math. Sci. \textbf{5} (2015) 13-18.




\bibitem[Sch69]{Schwartz}
J.~T.~Schwartz,
Nonlinear functional analysis,
Gordon and Breach science publishers, 1969.






\end{thebibliography}
\end{document}